\documentclass[12pt]{article}

\usepackage{amsmath,amssymb,enumerate,bbm}
\newcommand{\tmop}[1]{\ensuremath{\operatorname{#1}}}

\newtheorem{theorem}{Theorem}[section]

\newtheorem{lemma}[theorem]{Lemma}

\newtheorem{remark}[theorem]{Remark}

\numberwithin{equation}{section}
\newenvironment{proof}{\noindent\textbf{Proof\ }}{\hspace*{\fill}$\Box$\medskip}

 \makeatletter
 
 \newcommand{\Rmnum}[1]{\expandafter\@slowromancap\romannumeral #1@}
 \makeatother

\begin{document}

\title{On a Furstenberg-Katznelson-Weiss type theorem over finite
fields}\author{Le Anh Vinh\\
Mathematics Department\\
Harvard University\\
Cambridge, MA 02138, US\\
vinh@math.harvard.edu}\maketitle

\begin{abstract}
  Using Fourier analysis, Covert, Hart, Iosevich and Uriarte-Tuero (2008)
  showed that if the cardinality of a subset of the $2$-dimensional vector
  space over a finite field with $q$ elements is $\geqslant \rho q^2$, with
  $q^{- 1 / 2} \ll \rho \leqslant 1$ then it contains an isometric copy of
  $\geqslant c \rho q^3$ triangles. In this note, we give a graph theoretic
  proof of this result.
\end{abstract}

\section{Introduction}

A classical result due to Furstenberg, Katznelson and Weill (\cite{furstenberg}) says
that if $E \subset \mathbbm{R}^2$ has positive upper Lebesgue density, then
for any $\delta > 0$, the $\delta$-neighborhood of $E$ contains a congruent
copy of a sufficiently large dilate of every three point configuration. In \cite{covert},
Covert, Hart, Iosevich and Uriarte-Tuero investigated an analog of
this result in finite field geometries. They addressed the case of triangles
in two-dimensional vector spaces over finite fields.

Given $E \subset \mathbbm{F}_q^2$, where $\mathbbm{F}_q$ is a finite field of
$q$ elements, define
\begin{equation}
  T_3 (E) =\{(x, y, z) \in E \times E \times E\}/ \sim
\end{equation}
with the equivalence relation $\sim$ such that $(x, y, z) \sim (x', y', z')$
if there exists $\tau \in \mathbbm{F}_q^2$ and $O \in S O_2 (\mathbbm{F}_q)$,
the set of two-by-two orthogonal matrices over $\mathbbm{F}_q$ with
determinant $1$, such that
\begin{equation}
  (x', y', z') = (O (x) + \tau, O (y) + \tau, O (z) + \tau) .
\end{equation}
The main result of \cite{covert} is the following (see \cite{covert} and the references therein for the motivation and related results to this theorem).

\begin{theorem} (\cite{covert}) \label{old}
  Let $E \subset \mathbbm{F}_q^2$, and suppose that
  \begin{equation}
    |E| \geqslant \rho q^2
  \end{equation}
  for some $\frac{C}{\sqrt{q}} \leqslant \rho \leqslant 1$ with a sufficiently
  large constant $C > 0$. Then there exists $c > 0$ such that
  \begin{equation}
    |T_3 (E) | > c \rho q^3 .
  \end{equation}
\end{theorem}

In this note, we will provide a graph theoretic proof of this result. More
precisely, we only give a different proof of the key estimate (Estimate (\ref{t}) in Section 2) in Covert, Hart, Iosevich and Uriarte-Tuero's proof. Our result however is interesting in its
own right as it is related to the number of hinges (i.e. paths of length two) in a large subgraph of an
$(n, d, \lambda)$-graph. The rest of this note is organized as follows. In
Section 2, we study the arguments of Covert, Hart, Iosevich and Uriarte-Tuero
in \cite{covert} and discuss where graph theoretic methods can play a role. In
Section 3, we establish a theorem about the number of colored paths of length
two in a pseudo-random coloring of a graph. Using this theorem, we will give
another proof of Theorem \ref{old} in the last section.

\section{Covert, Hart, Iosevich and Uriarte-Tuero's arguments}

In this section, we follow closely the presentation in \cite{covert}. Covert, Hart,
Iosevich and Uriarte-Tuero observed that it suffices to show that if $|E| >
\rho q^2$, then
\begin{equation}\label{tp}
  |\{(a, b, c) \in \mathbbm{F}_q^3 : |T_{a, b, c} (E) | > 0\}| \geqslant c
  \rho q^3,
\end{equation}
where
\begin{equation}
  T_{a, b, c} (E) =\{(x, y, z) \in E \times E \times E : \|x - y\|= a, \|x -
  z\|= b, \|y - z\|= c\},
\end{equation}
with
\[ \|x\|= x_1^2 + x_2^2 . \]

This follows from the following lemma which states that over finite fields, a
(non-degenerate) simplex is defined uniquely (up to translation and rotation)
by the norms of its edges.

\begin{lemma} (cf. Lemma 2.1 in \cite{covert}) Let $P$ be a (non-degenerate) simplex with
  vertices $V_0, V_1, \ldots, V_k$ with $V_j \in \mathbbm{F}_q^d$. Let $P'$ be
  another (non-degenerate) simplex with vertices $V_0', \ldots, V_k'$. Suppose
  that
  \begin{equation}
    \|V_i - V_j \|=\|V_i' - V_j' \|
  \end{equation}
  for all $i, j$. Then there exists $\tau \in \mathbbm{F}_q^d$ and $O \in S
  O_d (\mathbbm{F}_q)$ such that $\tau + O (P) = P'$.
\end{lemma}

The key estimate of the proof of Theorem \ref{old} in \cite{covert} is the following result about hinges.

\begin{theorem} \label{old main} (cf. Theorem 2.2 in \cite{covert})
	Suppose that $E \subset \mathbbm{F}_q^2$ and
  let $a, b \neq 0$. Then
  \[ |\{(x, y, z) \in E \times E \times E : \|x - y\|= a, \|x - z\|= b\}| =
     |E|^3 q^{- 2} + O (q|E|) . \]
  And if $|E| \gg q^{3 / 2}$, then
  \begin{equation}\label{t}
   |\{(x, y, z) \in E \times E \times E : \|x - y\|= a, \|x - z\|= b\}| = (1
     + o (1)) |E|^3 q^{- 2} .
  \end{equation}   
\end{theorem}

Since $|E| \geqslant \rho q^2$ for some $\frac{1}{\sqrt{q}} \ll \rho \leqslant
1$, we have $|E| \gg q^{3 / 2}$. From (\ref{t}), by the pigeon-hole principle, there
exists $x \in E$ such that
\[ |\{(y, z) \in E \times E : \|x - y\|= a, \|x - z\|= b\}| \geqslant |E|^2
   q^{- 2} . \]
We have two cases.

Case 1. Suppose that the number of elements of $S O_2 (\mathbbm{F}_q)$ that
fix $x$ and keep $(y, z)$ inside the pinned hinge is no more than $\rho q$.
Since $|E| \geqslant \rho q^2$, the number of distinct distances $c$ from $\{y
\in E : \|x - y\|= a\}$ to $\{z \in E : \|x - z\|= b\}$ is at least
\[ |E|^2 q^{- 2} \frac{1}{\rho q} \geqslant \frac{1}{2} \rho q. \]
Since there are $(q - 1)^2$ possible choices for $a$ and $b$, (\ref{tp}) follows.

Case 2. Suppose that the number of elements of $S O_2 (\mathbbm{F}_q)$ that
fix $x$ and keep $(y, z)$ inside the pinned hinge is more than $\rho q$. Then
both the circle of radius $a$, centered at $x$, and the circle of radius $b$,
centered at $x$, contain more than $\rho q$ elements of $E$. It is shown that
(cf. Lemma 2.3 in \cite{covert}) the number of distinct distance $c$ from \ $\{y \in E : \|x -
y\|= a\}$ to $\{z \in E : \|x - z\|= b\}$ is at lest $\frac{\rho q}{4}$. Since
there are $(q - 1)^2$ possible choices for $a$ and $b$, (\ref{tp}) follows.

Thus the proof of Theorem \ref{old} can be reduced to Theorem \ref{old main}. The main purpose
of this note is to give a graph theoretic proof of Theorem \ref{old main}.

\section{Number of hinges in an $(n, d, \lambda)$-graph}

We call a graph $G = (V, E)$ $(n, d, \lambda)$-graph if $G$ is a $d$-regular
graph on $n$ vertices with the absolute values of each of its eigenvalues but
the largest one is at most $\lambda$. It is well-known that if $\lambda \ll d$
then an $(n, d, \lambda)$-graph behaves similarly as a random graph $G_{n, d /
n}$. Precisely, we have the following result (cf. Theorem 9.2.4 in \cite{alon-spencer}).

\begin{theorem} (\cite{alon-spencer}) \label{tool 1}
  Let $G$ be an $(n, d, \lambda)$-graph. For a vertex $v \in V$ and a subset
  $B$ of $V$ denote by $N (v)$ the set of all neighbors of $v$ in $G$, and
  let $N_B (v) = N (v) \cap B$ denote the set of all neighbors of $v$ in $B$.
  Then for every subset $B$ of $V$:
  \begin{equation}
    \sum_{v \in V} (|N_B (v) | - \frac{d}{n} |B|)^2 \leqslant
    \frac{\lambda^2}{n} |B| (n - |B|) .
  \end{equation}
\end{theorem}

The following result is an easy corollary of Theorem \ref{tool 1}

\begin{theorem} \label{tool 2}
  (cf. Corollary 9.2.5 in \cite{alon-spencer}) Let $G$ be an $(n, d, \lambda)$-graph. For every
  set of vertices $B$ and $C$ of $G$, we have
  \begin{equation}
    |e (B, C) - \frac{d}{n} |B\|C\| \leqslant \lambda \sqrt{|B\|C|},
  \end{equation}
  where $e (B, C)$ is the number of edges in the induced bipartite subgraph of
  $G$ on $(B, C)$ (i.e. the number of ordered pair $(u, v)$ where $u \in B$,
  $v \in C$ and $u v$ is an edge of $G$). 
\end{theorem}

Suppose that a graph $G$ of order $n$ is colored by $t$ colors. Let $G_i$ be
the induced subgraph of $G$ on the $i^{\tmop{th}}$ color. We call a
$t$-colored graph $G$ $(n, d, \lambda)$-r.c. (regularly colored) graph if $G_i$ is an $(n,
d, \lambda)$-graph for for each color $i \in \{1, \ldots, t\}$. The following
result gives us an estimate for the number of colored paths of length two in
an $(n,d,\lambda)$-r.c. graph $G$.

\begin{theorem}\label{main}
  Let $G$ be an $(n, d, \lambda)$-r.c. graph. For any two colors $r$, b and
  every set of vertices $E$ of $G$, we have
  \begin{equation}
    |e_{r, b} (E) - \left( \frac{d}{n} \right)^2 |E|^3 | \leqslant 2 \frac{\lambda d}{n} |E|^2 + \lambda^2|E|,
  \end{equation}
  where $e_{r, b} (E)$ is the number of $(r, b)$-colored paths of length two
  (i.e. the number of ordered triple $(u, v, w) \in E \times E \times E$ with
  $u v$, $v w$ are edges of $G$, $u v$ is colored $r$ and $v w$ is colored
  b). 
\end{theorem}

\begin{proof}
  For a vertex $v \in V$ let $N_E^r (v)$ and $N_E^b (v)$ denote the set of all
  $r$ neighbors and $b$ neighbors of $v$ in $E$, respectively. From Theorem
  \ref{tool 1}, we have
  \begin{eqnarray}
    \sum_{v \in E} (|N_E^r (v) | - \frac{d}{n} |E|)^2 \leqslant \sum_{v \in V}
    (|N_E^r (v) | - \frac{d}{n} |E|)^2 & \leqslant & \frac{\lambda^2}{n} |E|
    (n - |E|) \nonumber\\
    \sum_{v \in E} (|N_E^b (v) | - \frac{d}{n} |E|)^2 \leqslant \sum_{v \in V}
    (|N_E^b (v) | - \frac{d}{n} |E|)^2 & \leqslant & \frac{\lambda^2}{n} |E|
    (n - |E|) . \nonumber
  \end{eqnarray}
  Thus, by the Cauchy Schwarz inequality, we have
  \begin{eqnarray*}
    & & \left[ \sum_{v \in E} (|N_E^r (v) |  -  \frac{d}{n} |E|) ( \frac{d}{n} |E|  -
    |N_E^b (v) |) \right]^2\\ 
    & \leqslant & \left[ \sum_{v \in E} (|N_E^r (v) |
    - \frac{d}{n} |E|)^2 \right] \left[ \sum_{v \in E} (|N_E^b (v) | -
    \frac{d}{n} |E|)^2 \right] \leqslant \frac{\lambda^4}{n^2} |E|^2 (n - |E|)^2 .
  \end{eqnarray*}
  This implies that
  \begin{equation} \label{1}
    \left| \sum_{v \in E} N_E^r (v) N_E^b (v) + \left( \frac{d}{n} \right)^2
    |E|^3 - \frac{d}{n} |E| \sum_{v \in E} (N_E^r (v) + N_E^b (v)) \right|
    \leqslant \frac{\lambda^2}{n} |E| (n - |E|)
  \end{equation}
  From Theorem \ref{tool 2}, we have
  \begin{eqnarray}
    \left| \sum_{v \in E} N_E^r (v) - \frac{d}{n} |E|^2 \right| & \leqslant &
    \lambda |E| \label{2}\\
    \left| \sum_{v \in E} N_E^b (v) - \frac{d}{n} |E|^2 \right| & \leqslant &
    \lambda |E|. \label{3}
  \end{eqnarray}
  Putting (\ref{1}), (\ref{2}) and (\ref{3}) together, we have
  \begin{equation}
    \left| \sum_{v \in E} N^r_E (v) N_E^b (v) - \left( \frac{d}{n} \right)^2
    |E|^3 \right| \leqslant 2 \frac{\lambda d}{n} |E|^2 + \frac{\lambda^2}{n} |E| (n - |E|)
    < 2 \frac{\lambda d}{n} |E|^2 + \lambda^2|E|,
  \end{equation}
  completing the proof of the theorem.
\end{proof}

If we color an $(n,d,\lambda)$-graph by one color then Theorem \ref{main} becomes.

\begin{theorem} \label{path}
Let $G$ be an $(n, d, \lambda)$-graph. For every set of vertices $E$ of $G$, we have
  \begin{equation}
    |p_2(E) - \left( \frac{d}{n} \right)^2 |E|^3 | \leqslant 2 \frac{\lambda d}{n} |E|^2 + \lambda^2|E|,
  \end{equation}
  where $p_2(E)$ is the number of ordered paths of length two in $E$
  (i.e. the number of ordered triple $(u, v, w) \in E \times E \times E$ with
  $u v$, $v w$ are edges of $G$). In particular, if $|E| \gg \lambda \left( \frac{n}{d} \right)$ then the number of ordered paths of length two in $E$ is
    \begin{equation}
    (1 + o (1)) |E|^3 \left( \frac{d}{n} \right)^2 .
  \end{equation}
\end{theorem}

\begin{remark} 

Using the second moment method, it is not difficult to show that for every constant $p$ the
random graph $G (n, p)$ contains
\begin{equation}
  (1 + o (1)) p^r (1 - p)^{(^s_2) - r} \frac{n^s}{| \tmop{Aut} (H) |}
\end{equation}
induced copies of $H$. Alon extended this result to $(n, d, \lambda)$-graphs.
He proved that every large subset of the set of vertices of a $(n, d,
\lambda)$-graph contains the ``correct'' number of copies of any fixed small
subgraph (Theorem 4.10 in \cite{krivelevich-sudakov}).

\begin{theorem}\label{tool}
  (\cite{krivelevich-sudakov}) Let $H$ be a fixed graph with $r$ edges, $s$
  vertices and maximum degree $\Delta$, and let $G = (V, E)$ be an $(n, d,
  \lambda)$-graph, where, say, $d \leqslant 0.9 n$. Let $m < n$ satisfies $m
  \gg \lambda \left( \frac{n}{d} \right)^{\Delta}$. Then, for every subset $U
  \subset V$ of cardinality $m$, the number of (not necessarily induced) copies
  of $H$ in $U$ is
  \begin{equation}
    (1 + o (1)) \frac{m^s}{| \tmop{Aut} (H) |} \left( \frac{d}{n} \right)^r .
  \end{equation}
\end{theorem}

Note that, in the ``simple case'', $H$ is a path of length two, then Theorem \ref{tool} is weaker than Theorem \ref{path}. 
\end{remark} 

\section{Graph theoretic proof of (\ref{t})}

Let $\mathbbm{F}_q$ denote the finite field with $q$ elements where $q \gg 1$
is an odd prime power. For a fixed $a \in \mathbbm{F}_q^{\ast}$, the finite Euclidean graph $G_q(a)$ in $\mathbbm{F}_q^2$ is defined as the graph with vertex
set $\mathbbm{F}_q^2$ and the edge set
\[ E =\{(x, y) \in \mathbbm{F}_q^2 \times \mathbbm{F}_q^2 \mid x \neq y, ||x - y|| = a\}, \]
where $||.||$ is the analogue of Euclidean distance $||x|| = x_1^2+x_2^2$.
In \cite{medrano}, Medrano et al. studied the spectrum of these graphs and showed that these graphs are asymptotically Ramanujan graphs.  They proved the following result.

\begin{theorem} (\cite{medrano}) \label{tool 4}
The finite Euclidean graph $G_q(a)$ is regular of valency $q \pm 1$ for any $a \in \mathbbm{F}_q^{\ast}$. Let $\lambda$ be any eigenvalues of the graph $G_q(a)$ with $\lambda \neq$ valency of the graph then 
\begin{equation}
|\lambda| \leq 2\sqrt{q}.
\end{equation}
\end{theorem}

Now consider the set of colors $L =\{c_1, \ldots, c_{q - 1} \}$
corresponding to elements of $\mathbbm{F}_q^{\ast}$. We color the complete
graph $G_q = K_{q^2}$ with vertex set $\mathbbm{F}_q^2$ by $q - 1$ colors such that
$(x, y) \in \mathbbm{F}_q^2 \times \mathbbm{F}_q^2$ is colored by the color
$c_i$ if $||x - y|| = i$. Then from Theorem \ref{tool 4}, $G_q$ is a $(q^2,q \pm 1, 2\sqrt{q})$-r.c. graph. Estimate (\ref{t}) follows immediately from Theorem \ref{main}.  

\begin{remark}
Note that the conclusion of (\ref{t}) holds with the Euclidean norm $||.||$ is replaced by any non-degenerate quadratic form on $\mathbbm{F}_q^2$. This fact can be shown similarly as the above. Let $Q$ be a non-degenerate quadratic form on $\mathbbm{F}_q^2$. The finite Euclidean graph $E_q(Q,a)$ is defined as the graph with vertex set $\mathbbm{F}_q^2$ and the edge set
\begin{equation}
E = \{(x,y) \in \mathbbm{F}_q^2 \times \mathbbm{F}_q^2 \mid x \neq y, Q(x-y) = a \}.
\end{equation}

In \cite{bannai}, Bannai, Shimabukuro and Tanaka studied the spectrum of the graph $E_q(Q,a)$ and showed that these graphs are asymptotically Ramanujan graphs. 

\begin{theorem} (\cite{bannai}) \label{tool 5}
Let $Q$ be a non-degenerate quadratic form on $\mathbbm{F}_q^2$. The finite Euclidean graph $E_q(Q,a)$ is regular of valency $q \pm 1$ for any $a \in \mathbbm{F}_q^{\ast}$. Let $\lambda$ be any eigenvalues of the graph $E_q(Q,a)$ with $\lambda \neq$ valency of the graph then 
\begin{equation}
|\lambda| \leq 2\sqrt{q}.
\end{equation}
\end{theorem}

Similarly, consider the set of colors $L =\{c_1, \ldots, c_{q - 1} \}$
corresponding to elements of $\mathbbm{F}_q^{\ast}$. We color the complete
graph $G_q = K_{q^2}$ with vertex set $\mathbbm{F}_q^2$ by $q - 1$ colors such that
$(x, y) \in \mathbbm{F}_q^2 \times \mathbbm{F}_q^d$ is colored by the color
$c_i$ if $Q(x - y) = i$. Then from Theorem \ref{tool 4}, $G_q$ is a $(q^2,q \pm 1, 2\sqrt{q})$-r.c. graph. From Theorem \ref{tool 5} and Theorem \ref{main}, we are done. 

\end{remark}

\section*{Acknowledgments}

The research is performed during the author's visit at the Erwin Schr\"odinger International Institute for Mathematical Physics. The author would like to thank the ESI for hospitality and financial support during his visit.

\end{document}